\documentclass[preprint,12pt]{elsarticle}

\usepackage{amssymb}
\usepackage{amsthm,a4wide,tikz}
\usepackage{graphicx,epsfig}
\usepackage{amsmath,amsfonts,amsbsy}

\journal{}

\newtheorem{thm}{Theorem}
\newtheorem{cor}[thm]{Corollary}
\newtheorem{lem}[thm]{Lemma}

\newdefinition{rem}[thm]{Remark}
\newdefinition{defn}[thm]{Definition}
\newdefinition{exmp}[thm]{Example}
\newproof{pf}{Proof}

\begin{document}

\newcommand{\lp}{\left(}
\newcommand{\rp}{\right)}
\newcommand{\lsp}{\left[}
\newcommand{\rsp}{\right]}
\newcommand{\lop}{\left]}
\newcommand{\rop}{\right[}
\newcommand{\lbr}{\left\{}
\newcommand{\rbr}{\right\}}
\newcommand{\lang}{\langle}
\newcommand{\rang}{\rangle}

\newcommand{\A}{\mathcal{A}}
\newcommand{\Se}{\mathcal{S}}
\newcommand{\Q}{\mathcal{Q}}
\newcommand{\C}{\mathcal{C}}

\newcommand{\DA}{\mathcal{D}_{\A}}
\newcommand{\DS}{\mathcal{D}_{\Se}}
\newcommand{\DC}{\mathcal{D}_{\C}}

\newcommand{\AB}{\mathcal{A}_{B}^{\alpha,\beta}}
\newcommand{\SB}{\Se_{B}}
\newcommand{\QB}{\Q_{B}}
\newcommand{\CB}{\C_{B}}
\newcommand{\CEV}{\C_{EV}}

\newcommand{\DSB}{\mathcal{D}_{\SB}}
\newcommand{\DQB}{\mathcal{D}_{\QB}}
\newcommand{\DCB}{\mathcal{D}_{\CB}}

\newcommand{\Pic}{\mathcal{P}}

\newcommand{\Ext}{\textnormal{Ext}}

\newcommand{\Dex}{\Gamma}

\begin{frontmatter}



\title{\LARGE Extreme semilinear copulas}


\author[addr1]{Fabrizio Durante\corref{cor3}}
\ead{fabrizio.durante@unisalento.it}

\author[addr2]{Juan Fern\'andez-S\'anchez}
\ead{juanfernandez@ual.es}

\author[addr3]{Manuel \'Ubeda-Flores}
\ead{mubeda@ual.es}

\address[addr1]{Dipartimento di Scienze dell'Economia\\
Universit\`a del Salento, 73100 Lecce, Italy
}

\address[addr2]{
Research Group of Theory of Copulas and Applications\\
University of Almer\'ia, 04120 Almer\'ia, Spain}

\address[addr3]{
Department of Mathematics\\
University of Almer\'ia, 04120 Almer\'ia, Spain}

\cortext[cor3]{Corresponding author.}

\begin{abstract}
We study the extreme points (in the Krein-Milman sense) of the class of semilinear copulas and provide their characterization. Related results into the more general setting of conjunctive aggregation functions (i.e, semi--copulas and quasi--copulas) are also presented.
\end{abstract}


\end{frontmatter}

\section{Introduction}

In the study of complex systems it is of interest to synthesize the information coming from different sources into a single output, which is either numerical or represented by a suitable function, graph, etc. For instance, the copula representation has proved to be a suitable tool to describe uncertain inputs in a probabilistic framework (see, e.g., \cite{DurSem16,Nel06}) as well as in an imprecise setting (see, e.g., \cite{DurSpi10,Monetal15,OmlSto20}).

In order to represent various kinds of relationships among inputs, different families of copulas have been introduced and studied, mainly motivated by the question of identifying those copulas that may describe at the best some stylized facts of the problem at hand. Among these various families, we focus on the class of semilinear copulas, which have been introduced in \cite{Dur06JNPS,DKMS08FSS} and further investigated and generalized in \cite{DeBDeMMes07,DQUF07IS,FerUbe12,JwaDeBDeM16,JwaDeBDeM13,
MaiSchSch16} among others. Semilinear copulas can be constructed from their diagonal sections and, thus, their tail behaviour can be easily described (see, e.g., \cite{DurFerPap15}). Interestingly, this class has been characterized both from a probabilistic perspective, being the output of a stochastic models generated by different shocks (see also \cite{DurGirMaz16JCAM}), and from an analytical perspective, since its elements have sections that are linear on some specific segments of the unit square (see \cite{DKMS08FSS}).

The class of semilinear copulas is a convex and compact subset (under $L^{\infty}$ norm) of the class of copulas (see \cite{DKMS08FSS}) and, hence, by Krein--Milman Theorem \cite{AliBor06}, it is the closed convex hull of its extreme points.
We remind here that an extreme point $x$ of a convex set $A$ is a point $x\in A$ that is not an interior point of any line segment lying entirely in $A$. Thus, each element of a convex set $A$ can be approximated via linear combinations of the elements of ${\rm Ext}(A)$, i.e. the set of the extreme points of $A$.

Now, although the knowledge of extreme copulas can be of potential interest in the description of the whole class of copulas, even in the bivariate case, only a few examples of extreme copulas are available (e.g., shuffles of Min \cite{MikSheTay92,TruFer13JSPI}, hairpin copulas \cite{DurFerTru14}, extreme biconic copulas \cite{DurFerUbe19}), and a handle characterization of extreme copulas is still out of reach.

Our aim is, hence, to investigate the extreme elements in the class of semilinear copulas and to provide their characterization. Some consequences for the measurement of asymmetry maps in the class of semilinear copulas are also discussed. In order to complement our main results, we also consider this problem in the general setting of aggregation functions \cite{GMMP09}, by focusing on semi--copulas \cite{DurSem05Kyb} and quasi--copulas \cite{Arietal20,Sem17}.

\section{Preliminaries}

A (bivariate) {\it copula} is a distribution function concentrated on $[0,1]^2$ whose marginals are uniformly distributed on $[0,1]$. The importance of copulas in probability and statistics comes from {\it Sklar's theorem} \cite{Sklar1959}, which shows that the joint probability distribution $H$ of a pair of random variables and the corresponding marginal distributions $F$ and $G$ are linked by a copula $C$ in the following manner:
$$H(x,y)=C(F(x),G(y)),\quad\mbox{for every } x,y\in[-\infty,\infty].$$
If $F$ and $G$ are continuous, then the copula is unique; otherwise, the copula is uniquely determined on $\textnormal{Range}(F)\times\textnormal{Range}(G)$ (see, e.g., \cite{deAmo2012}). For a complete review on copulas, we refer to \cite{DurSem16,Nel06}.

A copula can be seen as a binary operation $C\colon[0,1]^{2}\longrightarrow [0,1]$ which satisfies:
\begin{enumerate}
\item the {\it boundary conditions} $C(t,0)=C(0,t)=0$ and $ C(t,1)=C(1,t)=t$ for every $t\in[0,1]$; and
\item the {\it 2-increasing property}, i.e., $V_C(R):=C(u_{2},v_{2})-C(u_{1},v_{2})-C(u_{2},v_{1})+C(u_{1},v_{1})\ge 0$, where $R=[u_1,u_2]\times [v_1,v_2]$ is a rectangle in $[0,1]^2$.
\end{enumerate}

We denote by ${\mathfrak{C}}$ the set of all bivariate copulas.


For any copula $C$ we have
$$W(x,y)=0\vee (x+y-1)\le C(x,y)\le x\wedge y=M(x,y),\quad \forall (x,y)\in[0,1]^2,$$
where $(a\vee b)=\max(a,b)$ and $(a\wedge b)=\min(a,b)$. The copulas $M$ and $W$ belong to Ext$({\mathfrak{C}})$ (or the set of extreme copulas), but $\Pi$, the copula of independent random variables---given by $\Pi(u,v)=uv$ for all $(u,v)\in[0,1]^2$---is not an extreme copula.

The {\it diagonal section} $\delta_{C}$ of a copula $C$ is the function defined by $\delta_{C}(t)=C(t,t)$ for every $t\in [0,1]$. It is characterized by the following conditions:
\begin{enumerate}
\item[(D1)] $\delta(0)=0$ and $\delta(1)=1$;
\item[(D2)] $\delta$ is non-decreasing;
\item[(D3)] $\delta(t)\le t$ for all $t\in [0,1]$; and
\item[(D4)] $|\delta(t')-\delta(t)|\le 2|t'-t|$ for all $t,t'\in [0,1]$.
\end{enumerate}
Any function that satisfies (D1)--(D4) is called \emph{diagonal} and the set of all diagonals is denoted by ${\mathfrak{D}}$. We recall that property (D4) is called $2$--Lipschitz condition and implies that a diagonal is absolutely continuous and almost everywhere (a.e) differentiable with respect to the Lebesgue measure $\lambda$.

\section{Extreme semilinear copulas}\label{sec:main}

A lower (respectively, upper) semilinear copula $C$ is an element of $\mathfrak{C}$ constructed from a linear interpolation between the values that $C$ assumes at the lower boundaries (respectively, upper boundaries) of the unit square and the values that $C$ assumes on the diagonal section (see \cite{DKMS08FSS}). Specifically, $C$ is called \emph{lower semilinear} if the mappings
\begin{align*}
& f_1\colon [0, x]\longrightarrow [0, 1],\quad f_1(t) := C(t, x),\\\
& f_2\colon [0, x]\longrightarrow [0, 1],\quad f_2(t) := C(x, t),
\end{align*}
are linear for all $x\in [0,1]$. As the \emph{survival copula} \cite{Nel06} of an upper semilinear copula $C$ is a lower semilinear copula (see, e.g., \cite{DKMS08FSS}), we will restrict our attention to lower semilinear copulas.

In closed form, lower semilinear copulas can be described in terms of their diagonal sections  by the expression
\begin{equation}\label{eq:semilinear}
C(u,v)=\frac{(u\wedge v)\,\delta_C(u\vee v)}{(u\vee v)}
\end{equation}
for all $(u,v)\in[0,1]^2$, with the convention $\frac{0}{0}:=0$.

Conversely, given $\delta\in \mathfrak{D}$, it may be of interest to characterize which conditions on $\delta$ ensure that a function of type \eqref{eq:semilinear} is a copula. This characterization is given in \cite[Theorem 4]{DKMS08FSS} and is recalled here.

\begin{thm}
The function given by \eqref{eq:semilinear} is a lower semilinear copula if, and only if, the functions $x\mapsto \varphi_\delta(x):=\delta(x)/x$ and $x\mapsto \eta_\delta(x):=\delta(x)/x^2$ are non-decreasing and non-increasing, respectively, on $]0,1]$.
\end{thm}

In the sequel, copulas of form \eqref{eq:semilinear} will be referred simply as semilinear copulas; its class will be denoted by ${\mathfrak{C}}_S$ (for a probabilistic interpretation of semilinear copulas, see \cite{Sloot20}). Moreover, ${\mathfrak{D}}_{{\mathfrak{C}}_S}$ will denote the set of diagonal sections of all the elements of ${\mathfrak{C}}_S$.

It is known that the set ${\mathfrak{D}}$ is convex and compact with respect to $L^{\infty}$ norm (see, e.g., \cite{DurFerPap15}). However, there is no simple relationship between extreme copulas --- i.e. extreme points of ${\mathfrak{C}}$ --- and extreme diagonals --- i.e. extreme points of ${\mathfrak{D}}$, except for the case $\delta(t)=t$ for all $t\in[0,1]$ for the copula $M$ (see \cite{DurFerUbe19}).

Interestingly, unlike the sets ${\mathfrak{C}}$ and ${\mathfrak{D}}$, we have a clear relationship between the set of extreme points of ${\mathfrak{C}}_{S}$, namely ${\rm Ext}(\mathfrak{C}_{S})$, and that of ${\mathfrak{D}}_{{\mathfrak{C}}_S}$, namely ${\rm Ext}(\mathfrak{D}_{{\mathfrak{C}}_{S}})$, as the following result shows.

\begin{thm}\label{T:main}
Let $\delta\in{\mathfrak{D}}_{{\mathfrak{C}}_S}$, and let $C_\delta$ be the corresponding semilinear copula given by \eqref{eq:semilinear}. Then, $C_\delta\in {\rm Ext}(\mathfrak{C}_{S})$ if, and only if, $\delta\in {\rm Ext}(\mathfrak{D}_{{\mathfrak{C}}_{S}})$.
\end{thm}
\begin{proof}
The proof is a direct consequence of the fact that
semilinear copulas keep convex combinations (see \cite[section 6]{DKMS08FSS}), i.e. $C_{\alpha\delta_1+(1-\alpha)\delta_2}=\alpha C_{\delta_1}+(1-\alpha) C_{\delta_2}$ for every $\alpha\in [0,1]$.
\end{proof}

To summarize, the sets $\mathfrak{C}_{S}$ and ${\mathfrak{D}}_{{\mathfrak{C}}_S}$ are compact and convex subsets, respectively, of $\mathfrak C$ and $\mathfrak D$,  both equipped with $L^\infty$--norm. The mapping
$$
\mathfrak{C}_{S}\to {\mathfrak{D}}_{{\mathfrak{C}}_S},\quad C\mapsto \delta_C,
$$
is a homeomorphism. Moreover, as a consequence of Theorem \ref{T:main}, to compute the extreme points of ${\mathfrak{C}}_{S}$, we only need to find the extreme points of ${\mathfrak{D}}_{{\mathfrak{C}}_S}$.

To this end, we present here some properties related to any diagonal $\delta\in{\mathfrak{D}}_{{\mathfrak{C}}_S}$. In the sequel, when we consider the derivative of a diagonal, we will refer to the points where it exists (we recall that such a derivative exists a.e.). Moreover, the inequalities in which the derivative appears must be understood almost everywhere.

We start by considering the following results that will be helpful in the sequel.
\begin{itemize}
\item Since $\delta \in {\mathfrak{D}}_{{\mathfrak{C}}_{S}}$ is absolutely continuous in $[t,s]$ with $0<t<s\leq 1$, then $x\mapsto\delta(x)/x$ and $x\mapsto\delta(x)/x^{2}$ are absolutely continuous on $[t,s]$. See, for instance, \cite[Theorem 7.1.10]{KK}.
\item Let $\delta \in {\mathfrak{D}}_{{\mathfrak{C}}_{S}}$. Then, since
$x\mapsto \delta(x)/x^2$ is non-increasing on \,$\lop 0,1\rsp$, we have $\delta(x)/x^2  \geq  \delta(1)/1^2=1$, which implies $\delta(x) >0$ on $]0,1[$. Moreover, since $\delta
(x)/x$ is non-decreasing, $x\delta ^{\prime }(x)-\delta (x)\geq 0$ for a.e. $%
x\in \,]0,1[$, which implies that $\delta ^{\prime }(x)>0$ for a.e. $x\in
\,]0,1[$. In addition, since $\delta (x)/x^{2}$ is non-increasing,
\begin{equation*}
x\delta ^{\prime }(x)-2\delta (x)\leq 0
\end{equation*}%
for a.e. $x\in \,]0,1[$, from which we deduce that $\delta ^{\prime }(x)\neq
2$ for a.e. $x\in \lbrack 0,1[$.
\end{itemize}

\bigskip
\begin{thm}\label{th:char_cop}
The extreme points of the set ${\mathfrak{D}}_{{\mathfrak{C}}_S}$ are the diagonal sections $\delta\in {\mathfrak{D}}_{{\mathfrak{C}}_S}$ such that
\begin{equation}\label{extremediagonals}
\lambda\left(\left\{x\in\,]0,1] \colon \delta'(x)\mbox{ exists and }\,\frac{1}{x}<\frac{\delta'(x)}{\delta(x)}<\frac{2}{x}\right\}\right)=0.
\end{equation}
\end{thm}
\begin{proof}
Let $\delta \in {\mathfrak{D}}_{{\mathfrak{C}}_{S}}$. We denote by $\Gamma$ the set of all $x\in \lbrack 0,1]$ such that $\delta ^{\prime }(x)$ exists.
Suppose
\begin{equation*}
\lambda \left( \left\{ {x\in }\Gamma \colon \frac{1}{x}<\frac{\delta
^{\prime }(x)}{\delta (x)}<\frac{2}{x}\right\} \right) >0.
\end{equation*}%
Then there exists $0<\varepsilon <1/2$ such that the set
\begin{equation*}
A\subseteq\left\{ {x\in }\Gamma \colon \frac{1}{x}+\varepsilon <\frac{\delta
^{\prime }(x)}{\delta (x)}<\frac{2}{x}-\varepsilon ,\,\delta ^{\prime
}(x)<2-\varepsilon ,\,\delta (x)<x-\varepsilon \right\}
\end{equation*}%
satisfies $a_{1}=\inf A>0$, $a_{2}=\sup A<1$, and $\lambda (A)>0$.

First, since $\delta(x)/x$ is non-decreasing, it follows that, $\delta (t)=t$ for some $t\in \,]0,1[$ implies $\delta(x)=x$ for $x\in [t,1]$. Thus,
\begin{equation*}
\left\{ {x\in }\Gamma \colon \frac{1}{x}<\frac{\delta ^{\prime }(x)}{\delta
(x)}<\frac{2}{x}\right\} \subseteq \left\{ {x\in \Gamma}\colon \,\delta (x)<x\right\}
\end{equation*}
On the other hand, since the set $\left\{x\in \Gamma \colon \,\delta^{\prime }(x)<2\right\}$ has Lebesgue measure $1$, it follows that
\begin{align*}
& \lambda \left( \left\{ {x\in }\Gamma \colon \frac{1}{x}<\frac{\delta
^{\prime }(x)}{\delta (x)}<\frac{2}{x},\,\delta ^{\prime }(x)<2,\,\delta
(x)<x\right\} \right)  \\
&=\lambda \left( \left\{ {x\in }\Gamma \colon \frac{1}{x}<\frac{\delta
^{\prime }(x)}{\delta (x)}<\frac{2}{x}\right\} \right) >0
\end{align*}

Since
$$
\left\{ {x\in }\Gamma \colon \frac{1}{x}<\frac{\delta
^{\prime }(x)}{\delta (x)}<\frac{2}{x},\,\delta ^{\prime }(x)<2,\,\delta
(x)<x\right\} =\bigcup_n A_{n}
$$
with
\begin{equation*}
A_{n}=\left\{ {x\in }\Gamma \colon \frac{1}{x}+\frac{1}{n}<\frac{\delta
^{\prime }(x)}{\delta (x)}<\frac{2}{x}-\frac{1}{n},\,\delta ^{\prime }(x)<2-%
\frac{1}{n},\,\delta (x)<x-\frac{1}{n}\right\} ,
\end{equation*}
it follows that there exists at least a set $A_{n_0}$ that satisfies $\lambda \left(A_{n_{0}}\right) >0$.

Finally, consider that it is possible to assume $a_{1}=\inf A>0$ and $a_{2}=\sup A<1$ by considering $A=A_{n_{0}}\cap \left[ \alpha ,1-\alpha
\right]$ with
\begin{equation*}
\alpha =\frac{\lambda \left( \left\{ {x\in }\Gamma \colon \frac{1}{x}<\frac{%
\delta ^{\prime }(x)}{\delta (x)}<\frac{2}{x}\right\} \right) }{4}.
\end{equation*}
Since $\lambda (A)>0,$ consider $A_{1}=A\cap \lbrack 0,b_{1}]$ and $A_{2}=A\cap \lbrack b_{2},1]$ such that
$$
\lambda \left( A_{1}\right) =\lambda \left(A_{2}\right) =\lambda(A) /4.
$$
Moreover, let $g$ be the function equal to $\mathbf{1}_{A_{2}}-\mathbf{1}_{A_{1}}$ -- here $\mathbf{1}_{B}$ is the characteristic function of the set $B$ -- and set $f:=\gamma g$, where $\gamma$ is a sufficiently small non-negative constant such that:
\begin{enumerate}
\item[(i)] $f$ is a measurable function, {$|f(x)|<\varepsilon /2$}, such that
$$
F(x):=\int_{0}^{x}f(t)\,dt\leq 0
$$
for every $x\in \lbrack 0,1]$, with $F(a_{1})=0$, {$%
F(a_{2})=0$}, and $-F(x)<\vartheta $, for $\mathbf{\vartheta =-}\ln\left( 1-\varepsilon /8\right)$;
\item[(ii)] if $x\in A$, then
\begin{equation*}
\frac{1}{x}<\frac{\delta ^{\prime }(x)}{\delta (x)}\pm \frac{f(x)}{%
2-e^{-\vartheta }}<\frac{2}{x};
\end{equation*}%
moreover, if $x\notin A$, then $f(x)=0$, and $\lambda \left( \left\{ x\in
\lbrack 0,1]:f(x)\neq 0\right\} \right) >0$;
\item[(iii)] $0<\delta ^{\prime }(x)\pm \delta (x)f(x)<2-\varepsilon /2$ for
all $x\in ]a_{1},a_{2}[$.
\end{enumerate}

Notice that the function $F\colon \lbrack 0,1]\rightarrow \mathbb{R}$ given above
is Lipschitz with constant $\gamma$.

Now, consider the function $\delta _{1}(x):=\delta (x)e^{F(x)}$. Thus, $%
\delta _{1}$ belongs to ${\mathfrak{D}}_{{\mathfrak{C}}_{S}}$. In fact, we
have:

\begin{itemize}
\item $\delta _{1}(0)=\delta (0)e^{F(0)}=0$ and $\delta _{1}(1)=\delta
(1)e^{F(1)}=1$; here observe that $F(x)=0$ for all $x\in \lbrack
0,a_{1}]\cup \lbrack a_{2},1]$.

\item $\delta _{1}(x)\leq xe^{F(x)}\leq x$ for all $x\in \lbrack 0,1]$.

\item $\delta _{1}^{\prime }(x)\geq 0$ for all $x\in \lbrack 0,1]$. In fact,
from (iii), we have $(\delta ^{\prime }(x)+\delta (x)f(x))e^{F(x)}\geq 0$
for $x\in ]a_{1},a_{2}[$; and $\delta _{1}(x)=\delta (x)$ for $x\in \lbrack
0,a_{1}]\cup \lbrack a_{2},1]$).

\item If $x\in \,]a_{1},a_{2}[$, then using (iii) we have almost everywhere that
\begin{equation*}
\delta _{1}^{\prime }(x)=(\delta ^{\prime }(x)+\delta (x)f(x))e^{F(x)}\leq
\delta ^{\prime }(x)+\delta (x)f(x)<2-\varepsilon /2<2;
\end{equation*}%
moreover, if $x\in \lbrack 0,a_{1}]\cup \lbrack a_{2},1]$, then {$\delta
_{1}^{\prime }(x)=\delta ^{\prime }(x)<2$.}
\end{itemize}

In order to prove that the function $\delta _{1}(x)/x$ is non-decreasing, we
consider the function $\ln \left( \delta _{1}(x)/x\right) $. First, we
observe that it is absolutely continuous for every interval of type $[t,1]$ with $0<t$ (see \cite{KK}) and, from (ii) we have almost everywhere that
\begin{equation*}
\left( \ln \left( \frac{\delta _{1}(x)}{x}\right) \right) ^{\prime }=\frac{%
\delta ^{\prime }(x)}{\delta (x)}+f(x)-\frac{1}{x}\geq 0.
\end{equation*}
A similar reasoning leads us to the fact that $\delta _{1}(x)/x^{2}$ is
non-increasing.

Second, consider the function $\delta _{2}(x):=\delta (x)\left(
2-e^{F(x)}\right) $. Thus, $\delta _{2}$ belongs to ${\mathfrak{D}}_{{%
\mathfrak{C}}_{S}}$. In fact, we have:

\begin{itemize}
\item $\delta _{2}(0)=\delta (0)\left( 2-e^{F(0)}\right) =0$ and $\delta
_{2}(1)=\delta (1)\left( 2-e^{F(1)}\right) =1$.

\item $\delta_2(x)=\delta(x)\le x$ for every $x\in [0,a_1]\cup[a_2,1]$, and,
for every $x\in]a_1,a_2[$, $\delta_2(x)<(x-\varepsilon)\left(2-e^{-%
\vartheta}\right)<x$ for a sufficiently small $\vartheta$.

\item $\delta_2^{\prime}(x)=\delta^{\prime}(x)\left(2-e^{F(x)}\right)-%
\delta(x)f(x)e^{F(x)}\ge \delta^{\prime}(x)-\delta(x)|f(x)|\ge 0$ for almost every $x\in]a_1,a_2[$ (from (iii)), and $\delta_2(x)=\delta(x)\ge 0$ for almost every $x\in [0,a_1]\cup[a_2,1]$.

\item For every $x\in]a_1,a_2[$, taking $e^{-\vartheta}=1-\varepsilon/8$, $|f(x)|<\varepsilon/2$, we have almost everywhere that
\begin{eqnarray*}
\delta_2^{\prime}(x)\!\!\!&=&\!\!\!\delta^{\prime}(x)\left(2-e^{F(x)}%
\right)-\delta(x)f(x)e^{F(x)}\le
(2-\varepsilon)\left(2-e^{-\vartheta}\right)+|f(x)| \\
&\le&\!\!\!(2-\varepsilon)\left(-1+\frac{\varepsilon}{8}\right)+\frac{%
\varepsilon}{2}=2-\frac{\varepsilon}{4}+\frac{\varepsilon^2}{8}<2;
\end{eqnarray*}
and, if $x\in[0,a_1]\cup[a_2,1]$, then $\delta_2^{\prime}(x)=\delta^{%
\prime}(x)< 2$.
\end{itemize}

In order to prove that the function $\delta _{2}(x)/x$ is non-decreasing,
observe that the function $\ln \left( \delta _{2}(x)/x\right) $ is absolutely continuous or every interval of type $[t,1]$ with $t>0$, and, from (ii) we have
\begin{equation*}
\left( \ln \left( \frac{\delta _{2}(x)}{x}\right) \right) ^{\prime }=\frac{%
\delta ^{\prime }(x)}{\delta (x)}-\frac{f(x)}{2-e^{F(x)}}-\frac{1}{x}\geq 0.
\end{equation*}%
A similar reasoning leads us to the fact that $\delta _{2}(x)/x^{2}$ is non-increasing. Summarizing, the above considerations implies that the extreme diagonals fulfil \eqref{extremediagonals}.

Now, let $\delta (x)=\left( \delta _{3}(x)+\delta _{4}(x)\right) /2$ with $\delta_3$ and $\delta_4$ belonging to $ {\mathfrak{D}}_{{\mathfrak{C}}_{S}}$. Since $
\delta _{i}^{\prime }(x)/\delta _{i}(x)\in \left[ 1/x,2/x\right] $ for $i=3,4$, then we have that
\begin{equation*}
\frac{\delta ^{\prime }(x)}{\delta (x)}=\frac{\displaystyle\frac{\delta
_{3}^{\prime }(x)}{2}+\displaystyle\frac{\delta _{4}^{\prime }(x)}{2}}{%
\displaystyle\frac{\delta _{3}(x)}{2}+\displaystyle\frac{\delta _{4}(x)}{2}}=%
\frac{\delta _{3}^{\prime }(x)+\delta _{4}^{\prime }(x)}{\delta
_{3}(x)+\delta _{4}(x)}=\frac{1}{x}.
\end{equation*}%
if, and only if, $\delta _{i}^{\prime }(x)/\delta _{i}(x)=1/x$ for $i=3,4$. In fact, the sufficient part is immediate. For the necessary part, if $\frac{\delta
_{3}^{\prime }(x)}{\delta _{3}(x)}<\frac{\delta _{4}^{\prime }(x)}{\delta
_{4}(x)},$ then
$$
\frac{\delta _{3}^{\prime }(x)}{\delta
_{3}(x)}<\frac{\delta _{3}^{\prime }(x)+\delta _{4}^{\prime }(x)}{\delta
_{3}(x)+\delta _{4}(x)}<\frac{\delta _{4}^{\prime }(x)}{\delta _{4}(x)}.
$$
Thus, it holds that $\frac{\delta _{3}^{\prime}(x)}{\delta _{3}(x)}<\frac{1}{x},$ which is absurd. It follows that $\frac{\delta _{3}^{\prime }(x)}{\delta _{3}(x)}=\frac{\delta _{4}^{\prime }(x)}{\delta _{4}(x)}=\frac{1}{x}$. Similarly, we have that $\delta ^{\prime }(x)/\delta (x)=2/x$ if, and only if, $\delta
_{i}^{\prime }(x)/\delta _{i}(x)=2/x$ for $i=3,4$. Since $\delta(x)>0$ when $x>0$, it holds that $\ln(\delta(x))$ is absolutely continuous in $[t,1]$ for $t>0$. In other words,
$$
\ln\left(\delta(t)\right)=-\int_{\left[ t,1\right] }\frac{\delta
^{\prime }(x)}{\delta (x)}dx.
$$
Therefore, $\delta (x)=\delta _{3}(x)=\delta_{4}(x)$, i.e. $\delta $ is an extreme point, and this completes the proof.
\end{proof}

\section{A subclass of extreme semilinear copulas}\label{Studydelta}

\begin{figure}[htb]
\begin{center}
\includegraphics[width=6cm]{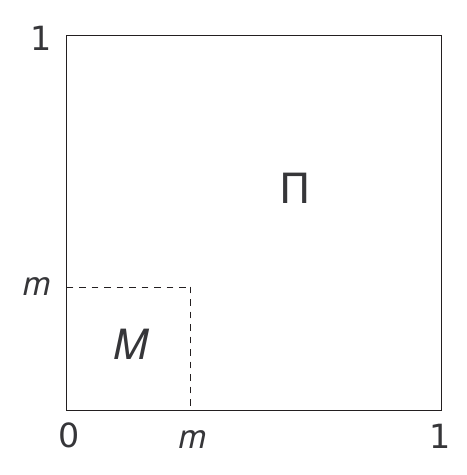}
\caption{Support of the extreme semilinear copulas generated by \eqref{eq:deltam}.}\label{fig:example}
\end{center}
\end{figure}

Here we present some examples of extreme semilinear copulas and study its closed convex hull. Specifically, we consider the elements of ${\rm Ext}(\mathfrak{C}_{S})$ that are generated by the following diagonal section
\begin{equation}\label{eq:deltam}
\delta_m(t)= (mt)\vee t^2\quad \mbox{for every }t\in[0,1],
\end{equation}
with $m\in[0,1]$. The support of the extreme semilinear copulas generated by \eqref{eq:deltam} can be obtained from Figure \ref{fig:example}.
These copulas can be interpreted in terms of rectangular patchwork construction, where the independence copula $\Pi$ is the background measure (see \cite{DeBDeM07IEEE,DurFerSem13IME,DurSamSar09CS}). Furthermore, such copulas include, as special cases, the independence copula $\Pi$ and the comonotonicity copula $M$.

Now, consider the closed convex hull of the set of diagonals $\delta_m$ given in eq. \eqref{eq:deltam}, denoted by ${\mathfrak{D}}_{\mathfrak{C}_{\delta_{m}}}$. We denote by $\mathfrak{C}_{\delta_{m}}$ the corresponding class of semilinear copulas generated by elements of ${\mathfrak{D}}_{\mathfrak{C}_{\delta_{m}}}$. To provide additional insights into the description of $\mathfrak{C}_{\delta_{m}}$, we remind the two following results from Functional Analysis (see, e.g., \cite{AliBor06}), which give a way to approximate and to represent elements of a compact convex set.

\bigskip
\begin{thm}[Krein-Milman]\label{T:KreinMilman}
Let $S$ be a non-empty compact convex subset of a locally convex Hausdorff topological vector space. Then $S$ is the closure of the convex hull of the set of extreme points of $S$.
\end{thm}

\begin{thm}[Choquet]\label{T:Choquet}
For a compact convex subset $K$ of a normed space $V$, given $k\in K$, there exists a probability measure $\nu$ supported on $\Ext(K)$ such that, for any affine function $f$ on $K$, we have
$$ f(k)=\int_K f(e)\, d\nu(e).$$
\end{thm}

As a consequence of Theorem \ref{T:KreinMilman}, the following result holds.

\begin{cor}
The set of semilinear copulas whose diagonal section $\delta$ is of type
\begin{equation}\label{eq:delta_dense}
\delta (x)=\left\{
\begin{array}{ll}
\alpha_{0}x, & 0\leq x\leq \alpha_{0}=x_{1}, \\
\alpha_{1}x+\beta_{1}x^{2}, & x_{1}\leq x\leq x_{2}, \\
\quad\vdots & \quad\vdots \\
\alpha_{j-1}x+\beta_{j-1}x^{2}, & x_{j-1}\leq x\leq x_{j}, \\
x^{2}, & x_{j}\leq x\leq 1,
\end{array}
\right.
\end{equation}
with $\alpha_k\in [0,1]$ for $k=0.\dots,j-1$ and $\beta_k\in [0,1]$ for $k=1,\dots,j-1$, is a dense subset in $\mathfrak{C}_{\delta_{m}}$.
\end{cor}
\begin{proof}
Because of the homeomorphism between $\mathfrak{C}_{\delta_{m}}$ and ${\mathfrak{D}}_{\mathfrak{C}_{\delta_{m}}}$, we only need to consider the problem in this latter set. First, we notice that the convex combinations of extreme points of ${\mathfrak{D}}_{\mathfrak{C}_{\delta_{m}}}$ can be represented in the form
\begin{equation*}
\delta^\ast(x)=\sum_{i=1}^{j}c_i\delta_{m_i}(x)=\left\{
\begin{array}{ll}
\biggl(\sum\limits_{i=1}^{j}c_im_i\biggr)x, & 0\leq x\leq m_1, \\
\biggl(\sum\limits_{i=2}^{j}c_im_i\biggr)x+c_1x^2, & m_{1}\leq x\leq m_{2}, \\
\quad\vdots & \quad\vdots \\
c_jm_jx+\biggl(\sum\limits_{i=1}^{j-1}c_i\biggr)x^2, & m_{j-1}\leq x\leq m_j, \\
x^{2}, & m_j\leq x\leq 1,
\end{array}
\right.
\end{equation*}
with $\sum_{i=1}^{j}c_i=1$. (Here, we have assumed, without loss of generality, that $i_1<i_2$ implies $m_{i_1}\le m_{i_2}$). The set of all such $\delta^\ast$ is dense in ${\mathfrak{D}}_{\mathfrak{C}_{\delta_{m}}}$ (since its closure coincides with the whole set because of Theorem \ref{T:KreinMilman}). Moreover, notice that each $\delta^\ast$ can be represented in the form \eqref{eq:delta_dense} by setting $\alpha_h=\sum_{i=h+1}^{j}c_im_i$ and $\beta_h=\sum_{i=1}^{h}c_i$, for $1\le h\le j$, and with $\alpha_j=0$. Now, the assertion follows by the fact that the homeomorphism $\mathfrak{C}_{\delta_{m}}\to {\mathfrak{D}}_{\mathfrak{C}_{\delta_{m}}}$ preserves convex combinations.
\end{proof}

Thus, elements of $\mathfrak{C}_{\delta_{m}}$ can be approximated via semilinear copulas generated by piecewise quadratic diagonal sections.

The next two results concerning semilinear copulas in $\mathfrak{C}_{\delta_{m}}$ are a consequence of Theorems \ref{T:main} and \ref{T:Choquet}.

\begin{cor}
The set $\mathfrak{C}_{\delta_{m}}$ is isomorphic to the set of all probability measures on ${\mathcal{B}}([0,1])$, the Borel subsets of $[0,1]$.
\end{cor}
\begin{proof}
Given $C\in\mathfrak{C}_{\delta_{m}}$, from Theorem \ref{T:Choquet} we have that
\begin{equation}\label{eq:semilinear-measure}
C(u,v)=\int_{0}^{1}C_{\delta_m}(u,v)\,d\mu(m),
\end{equation}
where $\mu$ is a probability measure in ${\mathcal{B}}([0,1])$, and the result follows.
\end{proof}

\begin{cor}
Let $\delta_m$ be the diagonal section given in \eqref{eq:deltam}. Then $\delta_m$ can be expressed in a unique way as
\begin{equation*}
\delta_m \left( t\right) =t^{2}F_{\mu }(t)+t\int_{t}^{1}m\,d\mu (m).
\end{equation*}
Moreover, the copula $C_{\delta_m}$ associated with $\delta_m$ can be written as
$$
C_{\delta_m}(u,v)=\left( \left( u\vee v\right) F_{\mu}(u\vee v)+\int_{u\vee v}^{1}m\,d\mu (m)\right) \left( u\wedge v\right),
$$
where $\mu $ is a probability measure in ${\mathcal{B}}([0,1])$, and $F_{\mu }$ is its distribution function.
\end{cor}

In view of the possible use in statistical applications, values of some popular measures of association for semilinear copulas have been considered in \cite{Dur06JNPS}. Here, we exploit the Choquet representation of a semilinear copula of type \eqref{eq:semilinear-measure} to obtained the desired results.

To this end, we consider three of the most common nonparametric measures of association between the components of a continuous random pair $(X,Y)$ are
{\it Spearman's rho} ($\rho$), {\it Gini's gamma} ($\gamma$), and the {\it Spearman's footrule coefficient} ($\varphi$). Such measures depend only on the
copula $C$ of $(X,Y)$, and are defined, e.g., in \cite{Nel06}.

The following result, whose proof is simple, provides the expressions for these measures when we consider an extreme semilinear copula with diagonal
$\delta_m$ given in eq. \eqref{eq:deltam}.

\begin{thm}\label{T:association}
Let $C_{\delta_m}$ be the semilinear copula with diagonal section given by \eqref{eq:deltam}. Then we have $\rho(C_{\delta_m})=m^{4}$,
$\varphi(C_{\delta_m})=m^{3}$, and
$$
\gamma(C_{\delta_m})=\begin{cases}
2m^3/3, &\qquad \mbox{if $m\le 1/2$}, \\
-2m^3/3+4m^2-3m+2/3, &\qquad \mbox{if $m>1/2$}.
\end{cases}
$$
\end{thm}

As a consequence of Theorems \ref{T:Choquet} and \ref{T:association}, we have the following result.

\begin{cor}
Let $C_\delta\in\mathfrak{C}_{\delta_{m}}$ be the semilinear copula of type \eqref{eq:semilinear-measure} with associated probability measure $\mu$. Then
\begin{eqnarray*}
\rho({C_{\delta }})\!\!\!&=&\!\!\!\int_{0}^{1}m^{4}\,d\mu (m),\\
\gamma({C_{\delta }})\!\!\!&=&\!\!\!\int_{0}^{1/2}\frac{2m^3}{3}\,d\mu (m)+\int_{1/2}^{1}\left(-\frac{2m^3}{3}+4m^2-3m+\frac{2}{3}\right)\,d\mu(m),\\
\varphi({C_{\delta }})\!\!\!&=&\!\!\!\int_{0}^{1}m^{3}\,d\mu (m).
\end{eqnarray*}
\end{cor}


\begin{rem}
We want to note that another interesting example of extreme semilinear copula is that generated by the diagonal function given by
\begin{equation}\label{deltasmaxmin}
\delta^{(p)}(t)=\frac{t^2}{p}\wedge t
\end{equation}
for all $t\in[0,1]$ with $p\in]0,1]$. Similar results for $\delta^{(p)}$ to those provided for $\delta_m$ in this section can be done analogously.
\end{rem}

\begin{rem}The previous observation gives us the opportunity to answer a natural question: can the results obtained in Section \ref{sec:main} be directly applied to upper semilinear copulas? In general, this is not possible since the extreme points of lower and upper semilinear copulas do not coincide.
For instance, consider $\delta_1$ be a diagonal of an upper semilinear copula. Then $\delta_2(t)=\delta_1(1-t)+2t-1$ for all
$t\in[0,1]$ is a diagonal of a lower semilinear copula. In the case that the extreme diagonals of the lower and upper semilinear copulas coincide, it would imply that $\delta_1=\delta^{(p)}$ is a diagonal of an upper semilinear copula, where $\delta^{(p)}$ is the diagonal given by \eqref{deltasmaxmin}. Thus,
we have
$$\delta_2(t)=\delta^{(p)}(1-t)+2t-1=t\wedge \left(\frac{t^2}{p}+\left(1-\frac{1}{p}\right)(2t-1)\right),$$
which is a diagonal for a lower semilinear copula, but $\delta_2(t)/t$ is non-increasing for $1-p<t<\sqrt{1-p}$.
\end{rem}

\section{Asymmetry maps of semilinear copulas}\label{sec:asymmetry}

The semilinear copula $C_\delta$ given by Eq. \eqref{eq:semilinear} is {\it exchangeable}, i.e. $C_\delta(u,v)=C_\delta(v,u)$ for all $(u,v)\in[0,1]^2$
(see also \cite{DeBDeMMes07}, where a method for constructing possibly asymmetric semilinear copulas with a given diagonal section is given). However,
other symmetry properties of copulas are of interest as well, as considered in \cite{FerUbe19}.

Here, we study three asymmetry maps for the set of the semilinear copulas; namely, we consider the asymmetry map with respect to the opposite diagonal
for the points $(u,v)$ and $(1-v,1-u)$ in $[0,1]^2$, and another one with respect to the points $(u,v)$ and $(1-u,1-v)$ in $[0,1]^2$, i.e. with respect
to the point $(1/2,1/2)$. Finally, we consider the mapping for measuring radial asymmetry.

\begin{rem}
By the term {\it asymmetry map} for the opposite diagonal we denote the function that makes the point $(u,v)$ correspond the maximum of the values of
$$
|C(u,v)-C(1-v,1-u)|
$$ for a copula $C$. Since semilinear copulas are exchangeable (i.e. $C(u,v)=C(v,u)$ for every $(u,v)\in [0,1]^2$), it suffices to study the case in which $u<v$. Since the symmetry with
respect to the opposite diagonal applies the triangle with vertices $(0,0), (1/2,1/2), (1,0)$ in the triangle $(1,0), (1,1), (1/2,1/2)$, it suffices to study the triangle $(0,0), (1/2,1/2), (1,0)$.
In the other two cases, the asymmetry maps are studied in similar manner.
\end{rem}

As it will be shown in the following, when studying these asymmetry maps for the set ${\mathfrak{C}}_S$, one will always consider the set
${\rm Ext}({\mathfrak{C}}_S)$. This fact will be a consequence of the following classical result, which is recalled here (see \cite{AliBor06}).

\begin{thm}[Bauer Maximum Principle]\label{T:Bauer}
If $J$ is a compact convex subset of a locally convex Hausdorff space, then every upper semicontinuous convex function on $J$ has a maximizer that is an extreme point.
\end{thm}

We start by considering asymmetry maps with respect to the opposite diagonal and to the point $(1/2,1/2)$. For every copula $C\in{\mathfrak{C}}$, consider the quantity
\begin{equation*}\label{chi}
\chi_C(u,v)=C(u\vee(1-v),v\vee(1-u))-C(u\wedge(1-v),v\wedge(1-u))
\end{equation*}
for $(u,v)\in[0,1]^2$. For a fixed $(u,v)\in[0,1]^2$, we wonder about the values of
\begin{equation}\label{maxmin}
\min_{C_\delta\in{\mathfrak{C}}_S}\chi_C(u,v)\quad{\rm and}\quad \max_{C_\delta\in{\mathfrak{C}}_S}\chi_C(u,v),
\end{equation}
which help to quantify the minimal and maximal asymmetry with respect to the opposite diagonal of the class of semilinear copula.

In order to study these values, it suffices to fix the point $(u,v)$ in the triangle $T_1$ of vertices $(0,0)$, $(1/2,1/2)$ and $(1,0)$, whence
\begin{equation*}\label{chi2}
\chi_{C_\delta}(u,v)=C_\delta(1-v,1-u)-C_\delta(u,v).
\end{equation*}
Thus, we have the following result:

\begin{thm}\label{T:oppositediag}
Let $C_\delta\in{\mathfrak{C}}_S$. Then we have
\begin{equation*}
L(u,v)\le \chi_{C_\delta}(u,v)\le U(u,v),
\end{equation*}
for every $(u,v)\in[0,1]^2$, where
\begin{equation}\label{L}
L(u,v)=\left\{
\begin{array}{ll}\vspace{0.2cm}
\displaystyle\frac{1-u-v}{1-(u\wedge v)}, & v\le 1-u, \\
\displaystyle\frac{u+v-1}{u\vee v}, & v>1-u,
\end{array}
\right.
\end{equation}
and
\begin{equation}\label{U}
U(u,v)=\left\{
\begin{array}{ll}\vspace{0.2cm}
(1-u-v)(1-(u\wedge v), & v\le 1-u, \\
(u+v-1)(u\vee v), & v>1-u,
\end{array}
\right.
\end{equation}
\end{thm}

\begin{proof}
Assume $C_{\delta}\in{\mathfrak{C}}_S$ and let $(u,v)\in T_1$, i.e. $v\le u\le 1-v$. If we define the function $\varphi_{\delta}(u):=\delta(u)/u$ for all $u\in]0,1]$, then it is clear that, for $u_1,u_2\in[0,1]$, $u_1<u_2$ implies
\begin{equation}\label{ecuaciondelta}
\frac{\varphi_{\delta}(u_2)u_1}{u_2}\le \varphi_{\delta}(u_1)\le \varphi_{\delta}(u_2).
\end{equation}
Moreover, $u\le\varphi_{\delta}(u)\le 1$ and
$$\frac{\varphi_{\delta}(1-v)u}{1-v}\le \varphi_{\delta}(u)\le\varphi_{\delta}(1-v).$$
All these chains of inequalities lead us to the following:
\begin{eqnarray*}
(1-u-v)(1-v)\!\!\!&\le&\!\!\! (1-u-v)\varphi_{\delta}(1-v)=(1-u)\varphi_{\delta}(1-v)-v\varphi_{\delta}(1-v)\\
&\le&\!\!\! (1-u)\varphi_{\delta}(1-v)-v\varphi_{\delta}(u)\le (1-u)\varphi_{\delta}(1-v)-\frac{\varphi_{\delta}(1-v)uv}{1-v}\\
&\le&\!\!\!\varphi_{\delta}(1-v)\left(1-u-\frac{uv}{1-v}\right)\le 1-u-\frac{uv}{1-v}.
\end{eqnarray*}
Since $C_\delta(1-u,1-v)-C_\delta(u,v)=(1-u)\varphi_{\delta}(1-v)-v\varphi_{\delta}(u)$, we obtain the bounds $L$ and $U$ in $T_1$. By using similar arguments, the rest of the proof follows.
\end{proof}

\begin{figure}[htb]
\begin{center}
\includegraphics[width=12cm]{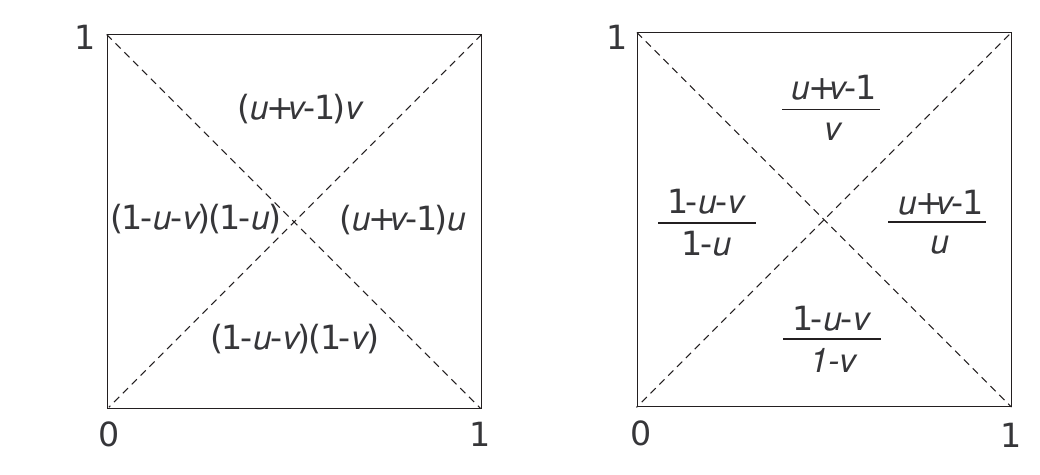}
\caption{The bounds $L$ (left) and $U$ (right) of Theorem \ref{T:oppositediag}.}
\end{center}\label{fig:bounds}
\end{figure}

The bounds $L$ and $U$ given by \eqref{L} and \eqref{U}, respectively, are shown in Figure \ref{fig:bounds}.

\begin{rem}Observe that, as a consequence of Theorem \ref{T:Bauer}, the minimum and maximum values in Theorem \ref{T:oppositediag} are reached for the
extreme diagonal sections of semilinear copulas given by
\begin{equation*}
\delta_{m}(t)=mt\vee t^2\quad{\rm and}\quad
\delta^{(p)}(t)=\frac{t^2}{p}\wedge t
\end{equation*}
for all $t\in[0,1]$, respectively, with $m\in[0,1]$ and $p\in]0,1]$.
\end{rem}

\begin{rem}
Notice that, although the maximum asymmetry with respect to the opposite diagonal is reached by other copulas, the copula $\Pi$ always gives us maximum asymmetry at every point $(u,v)\in[0,1]^2$.
\end{rem}

Now, for $C\in{\mathfrak{C}}$, consider the quantity given by
\begin{equation*}\label{xi}
\varrho_C(u,v)=C(u\vee (1-u), v\vee (1-v))-C(u\wedge (1-u),v\wedge (1-v))
\end{equation*}
for $(u,v)\in[0,1]^2$. This quantity measures the asymmetry of the copula $C$ with respect to the point $(1/2,1/2)$. Since every semilinear copula
$C_\delta$ is symmetric and  $\varrho_{C_\delta}(u,v)$ is linear for every $C_\delta\in{\mathfrak{C}}_S$, then we have similar results to those given
in Theorem \ref{T:oppositediag}.

\begin{cor}
Let $C_\delta\in{\mathfrak{C}}_S$. Then we have
\begin{equation*}
L(u,v)\le \varrho_{C_\delta}(u,v)\le U(u,v),
\end{equation*}
for every $(u,v)\in[0,1]^2$, where $L$ and $U$ are given by \eqref{L} and \eqref{U}, respectively.
\end{cor}

Finally, we consider a mapping for measuring radial asymmetry. We recall that a copula $C$ is {\it radially symmetric} if $C(u,v)=\widehat{C}(u,v)$ for all $(u,v)\in [0,1]^2$, where $\widehat{C}(u,v)=u+v-1+C(1-u,1-v)$ for every $(u,v)\in[0,1]^2$ is the {\it survival copula} associated with $C$ (see, e.g., \cite{Nel06}).

In \cite{Nel07}, it is proved that, for a given copula $C$ and every $(u,v)\in[0,1]^2$, a measure of radial asymmetry $\xi$ based on the $d_\infty$ distance can be defined as
$$\xi_C:=\underset{C\in \mathfrak{C}}{\max }\left\vert C(u,v)-\widehat{C}(u,v)\right\vert$$
i.e.
$$\xi_C:=\underset{C\in \mathfrak{C}}{\max }\left\vert C(u,v)-C(1-u,1-v)-u-v+1\right\vert$$
(see also \cite{Dehetal13,GenNes14,RosJoe13}).

We have the following result:

\begin{thm}\label{T:propradial}
Let $C_\delta\in{\mathfrak{C}}_S$. Then we have $\xi_{C_{\delta}}(u,v)\le U'(u,v)$, where
\begin{equation*}\label{Uprima}
U'(u,v)=\left\{
\begin{array}{ll}\vspace{0.2cm}
\displaystyle\frac{(1-u-v)(u\wedge v)}{1-(u\wedge v)}, & v\le 1-u, \\
\displaystyle\frac{(u+v-1)(1-(u\vee v))}{u\vee v}, & v>1-u.
\end{array}
\right.
\end{equation*}
\end{thm}

\begin{proof}
Assume $C_{\delta}\in{\mathfrak{C}}_S$. We define the function $\varphi_{\delta}(u):=\delta(u)/u$ for all $u\in]0,1]$ and consider two cases.
\begin{enumerate}
\item[(a)] $v\le u\le 1-v$. By using \eqref{ecuaciondelta} we have the following chain of inequalities:
\begin{eqnarray*}
-v(1-u-v)\!\!\!&\le&\!\!\! (1-u-v)(\varphi_{\delta}(1-v)-1)\\
&\le&\!\!\! (1-u)\varphi_{\delta}(1-v)-v\varphi_{\delta}(1-v)-1+u+v\\
&\le&\!\!\! (1-u)\varphi_{\delta}(1-v)-v\varphi_{\delta}(u)-1+u+v\\
&\le&\!\!\!(1-u)\varphi_{\delta}(1-v)-\frac{v\varphi_{\delta}(1-v)u}{1-v}-1+u+v\\
&\le&\!\!\! \left(1-u-\frac{uv}{1-v}\right)\varphi_{\delta}(1-v)-1+u+v\\
&\le&\!\!\! 1-u-\frac{uv}{1-v}-1+u+v=\frac{(1-u-v)v}{1-v}.
\end{eqnarray*}
Since $C_\delta(1-u,1-v)-C_\delta(u,v)-1+v+v=(1-u)\varphi_{\delta}(1-v)-v\varphi_{\delta}(u)-1+u+v$ then we have
\begin{eqnarray*}
|C_\delta(1-u,1-v)-C_\delta(u,v)-1+u+v|\!\!\!&=&\!\!\!\max\left(\frac{(1-u-v)v}{1-v},v(1-u-v)\right)\\
&=&\!\!\!\frac{(1-u-v)v}{1-v}.
\end{eqnarray*}
\item[(b)] $\max(v,1-v)\le u$. A similar reasoning to that in case (a) leads us to the following:
\begin{eqnarray*}
|C_\delta(1-u,1-v)\!\!\!&-&\!\!\!C_\delta(u,v)-1+u+v|\\
&=&\!\!\!\max\left(\frac{(1-u-v)(1-u)}{u},(1-u)(u+v-1)\right)\\
&=&\!\!\!\frac{(1-u-v)(1-u)}{u}.
\end{eqnarray*}
\end{enumerate}
By symmetry, the result easily follows.
\end{proof}

\begin{figure}
\begin{center}
\includegraphics[width=6cm]{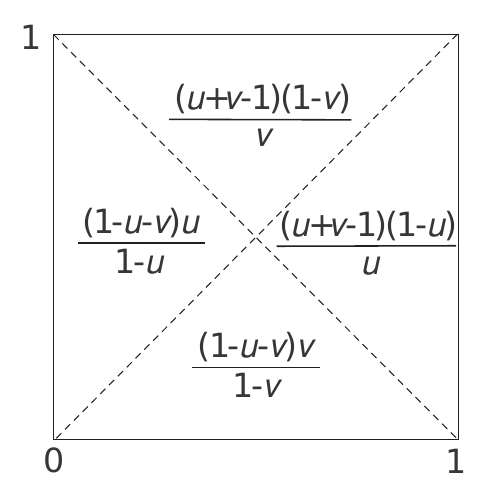}
\caption{The map of radial asymmetry of Theorem \ref{T:propradial}.}\label{fig:max_radial_asymm}
\end{center}
\end{figure}

The map of radial asymmetry map of Theorem \ref{T:propradial} is shown in Figure \ref{fig:max_radial_asymm}.

\begin{rem}Observe that, as a consequence of Theorem \ref{T:Bauer}, the maximum value in Theorem \ref{T:propradial} is reached for the extreme
(semilinear) diagonal $\delta^{(p)}$ given in Eq. \eqref{deltasmaxmin}.
\end{rem}

\section{Extreme semilinear semi-copulas and quasi-copulas}

The class of semilinear copulas can be extended to other two classes of aggregation functions, namely quasi--copulas and semi--copulas, as considered in \cite{Dur06JNPS, DKMS08FSS}.

We recall that a (bivariate) {\it semi-copula} is a function $S\colon [0,1]^2\longrightarrow [0,1]$ that is non-decreasing in each variable and admits uniform margins (see, e.g., \cite{DurQueSem06Kyb,DurSem05Kyb,DurSem16}); while a {\it quasi--copula} is a semi--copula that satisfies a Lipschitz property (see, e.g., \cite{Arietal20,GenQueRodSem99,Sem17}).


Semi-copulas of form \eqref{eq:semilinear} are characterized in terms of the properties of their diagonal sections by the following result (see \cite{DKMS08FSS}).

\begin{thm}
The function given by \eqref{eq:semilinear} is a semilinear semi-copula if, and only if, $\delta$ is non-decreasing, $0\le \delta(x)\le x$ for every $x\in [0,1]$, and $x\mapsto \delta(x)/x$ is non-decreasing on $]0,1]$.
\end{thm}

Here, we denote by ${\mathfrak{S}}_S$ the class of semilinear semi--copulas and by ${\mathfrak{D}}_{\mathfrak{S}_S}$ the class of their corresponding diagonal sections. It can be proved that ${\mathfrak{S}}_S$ is convex and compact in the topology of pointwise convergence (see, e.g., \cite{DurSem05Kyb}). Moreover, analogously to the copula case, Theorem \ref{T:main} also holds for the case of semilinear semi-copulas, so that the study of the set ${\rm Ext}({\mathfrak{S}}_S)$ is equivalent to study of the set ${\rm Ext}({\mathfrak{D}}_{\mathfrak{S}_S})$. Thus, we have the following result (whose proof is similar to that for biconic semi--copulas provided in \cite[Theorem 3.1]{DurFerUbe19} and, hence, it can be omitted here).

\begin{thm}
Let $\delta\in {\mathfrak{D}}_{\mathfrak{S}_S}$. Then $\delta\in Ext ({\mathfrak{D}}_{{\mathfrak{S}}_S})$ if, and only if, there exists $a\in [0,1]$ such that $\delta\in \left\{ \delta_{a}^R , \delta_{a}^L\right\}$, where
\begin{equation*}
\delta_{a}^R(x)=\left\{
\begin{array}{cc}
0, & 0\leq x<a, \\
x, & a\leq x\leq 1,
\end{array}
\right.
\quad and
\quad
\delta_{a}^L(x)=\left\{
\begin{array}{cc}
0, & 0\leq x\leq a, \\
x, & a<x\leq 1.
\end{array}
\right.
\end{equation*}
\end{thm}

In other words, diagonal sections of semilinear semi--copulas are only  the left- (respectively, right-) continuous step distribution functions with only one jump at $a\in[0,1]$.

\bigskip
Now, we consider the class of semilinear quasi--copulas (we denote this set by ${\mathfrak{Q}}_S$), whose characterization in terms of the respective diagonal sections is given by the following result (see  \cite{Dur06JNPS,DKMS08FSS}).

%

\begin{thm}\label{T:propdiagcuasi}
The function given by \eqref{eq:semilinear} is a semilinear quasi-copula if, and only if,
$\delta$ is a non-decreasing and $2$--Lipschitz function, $0\le \delta(x)\le x$ for every $x\in [0,1]$, and $\eta_\delta(x):=\delta(x)/x$ is non-decreasing in $]0,1]$ and satisfies
\begin{equation}\label{eq:diag-qc}
x_1\cdot\frac{\eta_\delta(x_2)-\eta_\delta(x_1)}{x_2-x_1}\le 1
\end{equation}
for every $x_1,x_2\in[0,1]$ with $x_1<x_2$.
\end{thm}

Notice that, condition \eqref{eq:diag-qc} is equivalent to
$\delta(x)\le x\delta'(x)\le x+\delta(x)$ at all points $x\in\, ]0,1]$ where the derivative $\delta'(x)$ exists (see also \cite[Corollary 18]{DKMS08FSS}).

We will denote by $\mathfrak{D}_{\mathfrak{Q}_S}$ the set of diagonal sections of semilinear quasi--copulas. We know that the set $\mathfrak{D}_{{\mathfrak{Q}}_S}$ is convex and compact with respect to $L^{\infty }$ norm.

In the sequel, when derivatives of diagonal sections in $\mathfrak{D}_{\mathfrak{Q}_S}$ are used, they are supposed to exist. In particular, the existence of such derivatives is guaranteed in a set of measure $1$, since these diagonals are absolutely continuous.

In order to study the set ${\rm Ext}({\mathfrak{Q}}_S)$, which is equivalent to study the set ${\rm Ext}({\mathfrak{D}}_{{\mathfrak{Q}}_S})$, we need a preliminary lemma.

\begin{lem}\label{L:lemacuasi}
Let $\delta\in {\mathfrak{D}}_{{\mathfrak{Q}}_S}$. Then we have $\delta (x)\geq x+x\ln(x)$ for all $x\in[0,1]$. Moreover, if there exists $x_0\in[0,1]$ such that $\delta(x_0)= x_0+x_0\ln(x_0)$, then $\delta(x)=x+x\ln(x)$ for all $x\in[x_0,1]$.
\end{lem}

\begin{proof}
First, note that the function $x\mapsto\delta (x)/x$ is absolutely continuous in $]0,1]$. From Theorem \ref{T:propdiagcuasi} we have
\begin{equation}\label{eq:diag-lemma}
\left( \frac{\delta (x)}{x}\right) ^{\prime }\leq \frac{1}{x}
\end{equation}
for all $x\in\,]0,1]$. Thus, for $x\in\,]0,1]$, we have
\begin{equation*}
\int_{x}^{1}\left( \frac{\delta (t)}{t}\right) ^{\prime }dt
\leq \int_{x}^{1}\frac{1}{t}\,dt,
\end{equation*}
whence we easily obtain $\delta (x)\geq x+x\ln(x)$.

On the other hand, if $x_1\in[x_0,1]$, then
\begin{equation*}
\frac{\delta(x_1)}{x_1}-\frac{\delta(x_0)}{x_0}=\int_{x_0}^{x_1}\left( \frac{\delta (t)}{t}\right) ^{\prime }dt
\leq \int_{x_0}^{x_1}\frac{1}{t}\,dt=\ln(x_1)-\ln(x_0).
\end{equation*}
Therefore, if $\delta(x_0)= x_0+x_0\ln(x_0)$, then it follows that $\delta(x_1)=x_1+x_1\ln(x_1)$, which completes the proof.
\end{proof}

Due to this fact, we can derive the following result.

\begin{thm}\label{T:char-qc}
The extreme points of the set $\mathfrak{D}_{{\mathfrak{Q}}_S}$ are the diagonal sections $\delta $ such that
\begin{equation}\label{eq:propextpoints}
\lambda \left( \left\{ x\in\, ] 0,1]: \left( \frac{\delta (x)}{x} \right) ^{\prime }=0\right\} \cup \left\{ x\in\, ] 0,1]: \left( \frac{\delta (x)}{x}\right) ^{\prime }=\frac{1}{x}\right\} \right) =1,
\end{equation}
where $\lambda$ is the Lebesgue measure on $[0,1]$.
\end{thm}

\begin{proof}Suppose $\delta \in {\mathfrak{D}}_{\mathfrak{Q}_{\mathfrak{S}}}$ be a diagonal satisfying \eqref{eq:propextpoints}. Suppose that there exist $\delta _{1},\delta _{2}\in D_{\mathfrak{Q}_{\mathfrak{S}}}$ such that $\delta(x) =\alpha \delta
_{1}(x)+(1-\alpha )\delta _{2}(x)$ for all $x\in[0,1]$ and for $\alpha\in ]0,1[$. We prove that $\delta_1=\delta_2=\delta$, from which we can conclude that $\delta \in{\rm Ext}( {\mathfrak{D}}_{\mathfrak{Q}_{\mathfrak{S}}})$. \\
It holds a.e.
$$
\left( \frac{\delta (x)}{x}\right)^{\prime }=\alpha \left( \frac{\delta _{1}(x)}{x}\right) ^{\prime }+(1-\alpha )\left( \frac{\delta _{2}(x)}{x}\right) ^{\prime }.
$$
Assume $\delta _{1}\neq \delta$ (the case $\delta _{2}\neq \delta$ is similar and we omit it) with
$$\lambda \left( \left\{ x\in]0,1]:\left( \frac{\delta
(x)}{x}\right) ^{\prime }\neq \left( \frac{\delta _{1}(x)}{x}\right)^{\prime }\right\} \right) >0.$$
We consider two cases.
\begin{itemize}
\item[(a)] $\displaystyle\lambda \left( \left\{x\in[0,1]:0=\left( \frac{\delta (x)}{x}\right) ^{\prime }\text{ and}\,\left( \frac{\delta _{1}(x)}{x}\right) ^{\prime }>0\right\} \right) >0.$\\
In this case, we have
    $$0=\alpha \left( \frac{\delta _{1}(x)}{x}\right)^{\prime }+(1-\alpha )\left( \frac{\delta _{2}(x)}{x}\right) ^{\prime }\quad {\rm a.e.}$$
    This implies
    $$\left( \frac{\delta _{2}(x)}{x}\right) ^{\prime}<0$$
    in a set of positive measure, but this contradicts the fact that $\delta _{2}(x)/x$ is non-decreasing.

\item[(b)] $\displaystyle\lambda \left( \left\{ x\in]0,1]:\left( \frac{\delta(x)}{x}\right) ^{\prime }=\frac{1}{x}\,\text{ and }\left( \frac{\delta _{1}(x)}{x
}\right) ^{\prime }<\frac{1}{x}\right\} \right) >0$.\\
In this case, we have
$$\frac{1}{x}=\alpha \left( \frac{\delta _{1}(x)}{x}\right) ^{\prime
}+(1-\alpha )\left( \frac{\delta _{2}(x)}{x}\right) ^{\prime }\quad {\rm a.e.}$$
This implies
$$\left( \frac{\delta _{2}(x)}{x}\right) ^{\prime }>\frac{1}{x}$$
in a set of measure positive, which contradicts the fact that
$$
\left( \frac{\delta _{2}(x)}{x}\right) ^{\prime }\leq \frac{1}{x}
$$
(recall Eq. \eqref{eq:diag-lemma}).
\end{itemize}
Therefore, we have
$$\left( \frac{\delta(x)}{x}\right) ^{\prime }=\left( \frac{\delta _{1}(x)}{x}\right) ^{\prime }=\left( \frac{\delta _{2}(x)}{x}\right) ^{\prime }$$
a.e. Since $\frac{\delta(x)}{x}$ and $\frac{\delta_1(x)}{x}$ are absolutely continuous on $[x,1]$, then we obtain
$$\frac{\delta(1)}{1}-\frac{\delta(x)}{x}=\int_{x}^{1}\left( \frac{\delta(t)}{t}\right) ^{\prime }dt=\frac{\delta_1(1)}{1}-\frac{\delta_1(x)}{x}.$$
We conclude that  $\delta _{1}=\delta$, whence $\delta $ is an extreme point.

Conversely, ab absurdo, suppose that $\delta \in{\rm Ext}( D_{\mathfrak{Q}_{\mathfrak{S}}})$ is a diagonal section which does not satisfy the condition in \eqref{eq:propextpoints}, so that
\begin{equation}\label{eq:diag-cond2}
\lambda \left( \left\{ x\in ]0,1]\colon 0<\left( \frac{\delta (x)}{x}\right) ^{\prime }<\frac{1}{x}\right\} \right) =\gamma >0.
\end{equation}
Suppose that there exists an interval $]0,r]$ in which $\delta(x)/x=a\ge 0$ for $x\in]0,r]$, and let $r_0$ be the supremum of the values $r$ with this property. If there is no interval $]0,r]$ in which $\delta(x)/x=a$ then $r_0=0$. Let $a_r$ be the solution of the equation $a=1+\ln(x)$. From Lemma \ref{L:lemacuasi} we have $1+\ln (r_0)\le a_r$ and, since the function $1+\ln(x)$ is increasing, we obtain $r_{0}\le a_r$. Moreover, we have $r_0<a_r$; otherwise, if $r_{0}=a_{r}$, then we would have $\delta (a_{r})/a_{r}=a=1+\ln(a_{r})$, i.e. $\delta (a_{r})=a_{r}+a_{r}\ln(a_{r})$. However, from Lemma \ref{L:lemacuasi} we have $\delta (x)=x+x\ln(x)$ for $x\in[a_{r},1]$, and hence $\gamma =0$, which is a contradiction.

In the case $r_{0}=1$, $\delta (x)/x=1$ for all $x\in[0,1]$, so that $\delta (x)=x$, which contradicts \eqref{eq:diag-cond2}.

In the case $r_{0}<1$, consider the value $r_{1}\in[0,1]$ such that either $\delta (x)=x$ or $\delta(x)=x+x\ln(x)$ for all $x\in\left[ r_{1},1\right] $, and $\delta(x)\neq x+x\ln(x)$ and $\delta(x)\neq x$ for $x\in[r_{1}-\varepsilon ,1]$ for a value $\varepsilon>0$. Then $r_{0}<r_{1}$. The next step is to modify
$$\left( \frac{\delta (x)}{x}\right) ^{\prime }$$
in the set
$$[r_{0}+\beta ,r_{1}-\beta ]\cap \left\{ x\in \lbrack 0,1]\colon 0<\left(\frac{\delta (x)}{x}\right) ^{\prime }<\frac{1}{x}\right\},$$
where $\beta \ge 0$ is such that
\begin{equation*}
\lambda \left( \lbrack r_{0}+\beta ,r_{1}-\beta ]\cap \left\{ x\in \lbrack 0,1]\colon 0<\left( \frac{\delta (x)}{x}\right) ^{\prime }<\frac{1}{x}\right\} \right) >0
\end{equation*}
(note that this is possible since $\gamma >0$). For $x\in[r_{0}+\beta ,r_{1}-\beta ]$ we have $\psi (x)<\delta (x)<x$, where $\psi (x):=0\wedge \left( x+x\ln(x)\right) $ for all $x\in ]0,1]$. Let
$$
\sigma :=\inf \left\{ x-\delta (x),\delta (x)-\psi (x):x\in\lbrack r_{0}+\beta ,r_{1}-\beta ]\right\}.
$$
Since the function $f(x):=x-\delta(x)$ is continuous, it reaches its minimum $x_0$ in the interval $\lbrack r_{0}+\beta ,r_{1}-\beta ]$; moreover, $x_0>0$ since $f(x)> 0$ for all $x\in \lbrack r_{0}+\beta ,r_{1}-\beta ]$. The same happens for the function $\delta (x)-\psi (x)$. Therefore, we have $\sigma>0$.

We define the functions $\delta _{1}(x)=\delta (x)+\Phi (x)$ and $\delta_{2}(x)=\delta (x)-\Phi (x)$ for every $x\in[0,1]$, where $\Phi (x)/x$ is an absolutely continuous function satisfying the following conditions:
\begin{enumerate}
\item[(i)] $\displaystyle\left( \frac{\Phi (x)}{x}\right) ^{\prime}=0$ provided that $x\in \lbrack 0,1]\backslash \displaystyle\left\{ x\in \lbrack 0,1]\colon 0<\left( \frac{\delta (x)}{x}\right) ^{\prime }<\frac{1}{x}\right\} ,$

\item[(ii)] $2\displaystyle\left\vert \left( \frac{\Phi (x)}{x}\right) ^{\prime }\right\vert \leq \min\left\{ \left( \frac{\delta (x)}{x}\right) ^{\prime },\frac{1}{x}-\left(\frac{\delta (x)}{x}\right) ^{\prime }\right\}$ a.e.,

\item[(iii)] $\Phi (x)=0$ provided that $x\in \lbrack 0,r_{0}+\beta ]\cup \lbrack r_{1}-\beta ,1],$

\item [(iv)] $\displaystyle\underset{x\in \lbrack 0,1]}{\max }\Phi (x)<\frac{\sigma}{2}.$
\end{enumerate}
To guarantee the existence of the function $\Phi$ satisfying the above conditions, observe that
\begin{align*}
\big\lbrack r_{0}+\beta,r_{1}-\beta \big]&\bigcap \displaystyle\left\{ x\in \lbrack 0,1]\colon 0<\left( \frac{\delta (x)}{x}\right) ^{\prime }<\frac{1}{x}\right\}\\
&=\bigcup_n\left(\lbrack r_{0}+\beta, r_{1}-\beta ]\cap \displaystyle\left\{ x\in \lbrack 0,1]\colon \frac{1}{n}<\left( \frac{\delta (x)}{x}\right) ^{\prime }<\frac{1}{x}-\frac{1}{n}\right\}\right).
\end{align*}
Thus, there exists $n_0\in\mathbb{N}$ such that, if
$$A:=\lbrack r_{0}+\beta, r_{1}-\beta ]\cap \displaystyle\left\{ x\in \lbrack 0,1]\colon \frac{1}{n_0}<\left( \frac{\delta (x)}{x}\right) ^{\prime }<\frac{1}{x}-\frac{1}{n_0}\right\},$$
then $\lambda(A)>0$. We divide $A$ into two measurable sets, $A_1$ and $A_2$, such that $\lambda(A_1)=\lambda(A_2)$ and define the function
$$\Phi(x):=\frac{x\sigma}{8 n_0}\int_{0}^{x}\left({\bf 1}_{A_{1}}(t)-{\bf 1}_{A_{2}}(t)\right)dt.$$
Then the function $\Phi$ satisfies conditions (i)--(iv).

Now, we check that $\delta _{1}\in {\mathfrak{D}}_{\mathfrak{Q}_{\mathfrak{S}}}$ (the proof for $\delta_2$ is similar and we omit it). From condition (iv) we have $0\leq \psi (t)\leq \delta (t)+t\Phi (t)\leq t$ for every $t\in\lbrack r_{0}+\beta ,r_{1}-\beta ]$.

Finally, we check
$$\frac{\delta _{1}(x)}{x}\leq \delta_{1}^{\prime }(x)\leq 1+\frac{\delta _{1}(x)}{x},$$
or equivalently,
$$0\leq \left( \frac{\delta _{1}(x)}{x}\right) ^{\prime }\leq\frac{1}{x},$$
but this follows from condition (ii) and the fact that
$$\frac{\delta (x)}{x}\leq \delta ^{\prime }(x)\leq 1+\frac{\delta (x)}{x}.$$
It is clear that $\delta(x) =(\delta _{1}(x)+\delta _{2}(x))/2$, i.e. $\delta $ is not an extreme diagonal, and this completes the proof.
\end{proof}

Notice that, if $\delta(t)=t^2$ on $[0,1]$, then $\delta$ is an extreme diagonal in $\mathfrak{D}_{{\mathfrak{C}}_S}$, but it is not an extreme diagonal in $\mathfrak{D}_{{\mathfrak{Q}}_S}$ since it does not satisfy the conditions of Theorem \ref{T:char-qc}.


The next example provides a family of diagonal sections of extreme quasi--copulas, which are not copulas.

\begin{exmp}
Consider the family of diagonals
\begin{equation*}
\delta _{\beta}(x)=\begin{cases}
\beta x, &\quad\mbox{if } 0\leq x<e^{\beta-1}, \\
x+x\ln(x), &\quad\mbox{if }  e^{\beta-1}\leq x\leq 1,
\end{cases}
\end{equation*}
where $\beta\in[0,1[$. It is easy to check that this family of diagonal sections satisfies \eqref{eq:propextpoints} and, therefore, it belongs to the set ${\rm Ext}({\mathfrak{D}}_{\mathfrak{Q}_{\mathfrak{S}}})$. Furthermore, none of the associated quasi-copulas $Q_\beta$ is a copula, since, for instance, $Q_\beta(u,u)=\beta u<u^2$ for every $u\in\,]\beta,e^{\beta-1}[$ (in fact, it is known from  \cite{Dur06JNPS,DKMS08FSS} that every semilinear copula $C$ satisfies $C(u,v)\ge uv$ for all $(u,v)\in[0,1]^2$). We also want to observe that, after some elementary calculations, it is easy to check that the semilinear quasi--copula associated with $\delta_{\beta}$ spread a negative mass equal to $\beta-1$  on the segment joining the points $\left( e^{\beta-1},e^{\beta-1}\right)$ and $\left( 1,1\right) $.
\end{exmp}

\section{Conclusions}

We have studied the extreme points of semilinear semi--copulas, quasi--copulas and copulas. In particular, we have proved that an extreme semilinear (semi--, quasi--)copula is characterized by the corresponding extreme diagonal section.

\section*{Acknowledgements}

The first author has been supported by the project ``Stochastic Models for Complex Systems'' by Italian MIUR (PRIN 2017, Project no. 2017JFFHSH).



\end{document}